\documentclass[12pt]{article}
\usepackage{epic,latexsym,amssymb}
\usepackage{amsfonts}
\usepackage{amscd}
\usepackage{amsmath}
\usepackage{graphicx}
\usepackage{color}
\usepackage{caption,subcaption}
\usepackage{tikz}

\usepackage{float}

\newtheorem{thm}{Theorem}
\newtheorem{conj}{Conjecture}
\newtheorem{ob}[thm]{Observation}

\newtheorem{quest}{Question}

\newtheorem{defn}{Definition}
\newtheorem{claim}{Claim}
\newtheorem{cor}[thm]{Corollary}

\newcommand{\pn}{{\rm pn}}

\newcommand{\proof}{\noindent\textbf{Proof. }}

\newcommand{\smallqed}{{\tiny ($\Box$)}}
\newcommand{\qed}{$\Box$}

\newcommand{\1}{\vspace{0.1cm}}

\newcommand{\cart}{\, \Box \,}
\newcommand{\cartr}{\, \Box_R \,}

\newcommand{\QEDmark}{\mbox{\textsc{qed}}}
\newcommand{\proofStarter}[1]{\textsc{#1} }

\def\vertex(#1){\put(#1){\circle*{2}}}
\def\vertexo(#1){\put(#1){\circle{2}}}
\def\vert(#1){\put(#1){\circle*{1.5}}}
\def\verto(#1){\put(#1){\circle{1.5}}}
\def\lab(#1)#2{\put(#1){\makebox(0,0)[c]{#2}}}
\definecolor{DarkGreen}{rgb}{0.2, 0.6, 0.3}

\definecolor{electricindigo}{rgb}{0.44, 0.0, 1.0}

\begin{document}

\title{On a Vizing-type integer domination conjecture}
\author{$^{1,2}$Randy Davila and $^3$Elliot Krop\\
\\
$^1$Department of Pure and Applied Mathematics\\
University of Johannesburg \\
Auckland Park 2006, South Africa \\
%\small {\tt Email: davilar@uhd.edu} \\
\\
$^2$Department of Mathematics and Statistics\\
University of Houston--Downtown \\
Houston, TX 77002 \\
\small {\tt Email: davilar@uhd.edu} \\
\\
$^3$Department of Mathematics \\
Clayton State University \\
Morrow, GA 30260, USA \\
\small {\tt Email: elliotkrop@clayton.edu}
}

\date{}
\maketitle

\begin{abstract}
Given a simple graph $G$, a dominating set in $G$ is a set of vertices $S$ such that every vertex not in $S$ has a neighbor in $S$. Denote the domination number, which is the size of any minimum dominating set of $G$, by $\gamma(G)$. For any integer $k\ge 1$, a function $f : V (G) \rightarrow \{0, 1, . . ., k\}$ is called a \emph{$\{k\}$-dominating function} if the sum of its function values over any closed neighborhood is at least $k$. The weight of a $\{k\}$-dominating function is the sum of its values over all the vertices. The $\{k\}$-domination number of $G$, $\gamma_{\{k\}}(G)$, is defined to be the minimum weight taken over all $\{k\}$-domination functions. 
Bre\v{s}ar, Henning, and Klav\v{z}ar (On integer domination in graphs and Vizing-like problems. \textit{Taiwanese J. Math.} \textbf{10(5)} (2006) pp. 1317--1328) asked whether there exists an integer $k\ge 2$ so that $\gamma_{\{k\}}(G\cart H)\ge \gamma(G)\gamma(H)$.
In this note we use the Roman $\{2\}$-domination number, $\gamma_{R2}$ of Chellali, Haynes, Hedetniemi, and McRae, (Roman $\{2\}$-domination. \textit{Discrete Applied Mathematics} \textbf{204} (2016) pp. 22-28.) to prove that if $G$ is a claw-free graph and $H$ is an arbitrary graph, then $\gamma_{\{2\}}(G\cart H)\ge \gamma_{R2}(G\cart H)\ge \gamma(G)\gamma(H)$, which also implies the conjecture for all $k\ge 2$.
\end{abstract}


\section{Introduction}
Given a graph $G$, a \emph{dominating set} is a set of vertices $S$ in $G$ with the property that every vertex not in $S$ has a neighbor in $S$. The \emph{domination number} of $G$, written $\gamma(G)$, is the cardinality of a minimum dominating set in $G$. One of the most influential and widely studied conjectures of domination in graphs is \emph{Vizing's conjecture}, originally posed by V.G. Vizing in 1963 \cite{Vizing}. This conjecture states that the domination number of the Cartesian product of graphs $G$ and $H$ is bounded below by the product of the domination numbers of $G$ and $H$. 
\begin{conj}{\rm (Vizing's Conjecture \cite{Vizing})}\label{Vizing}
For every pair of finite graphs $G$ and $H$,
\begin{equation*}
\gamma(G\cart H) \geq \gamma(G)\gamma(H).
\end{equation*}
\end{conj}

Attemps to prove Vizing's conjecture are numerous. At this point, the conjecture is known for many pairs of graphs, but the strongest results of this type are those where one of the graphs, say $G$, is taken from a fixed family and the other graph $H$ is chosen arbitrarily. A broad result of this kind was that of Bartsalkin and German from 1979. They showed that the conjecture holds when $G$ is a ``decomposable'' graph (we won't list the definition here), as well as a spanning subgraph of decomposable graphs which shares the same domination number as the parent graph. Graphs that satisfy this defintion are said to belong to \emph{class $A_0$}. It is not clear which classes of graphs are included in class $A_0$ but it is not hard to show that all trees are as well as graphs with two packing number equal to their domination number. A slight extention of class $A_0$ was made by Hartnell and Rall in 1995 and is known as type $\chi$. In 2009, Aharoni and Szabo proved Vizing's conjecture for chordal graphs using previous work on the independence-domination number from matching theory. The conjecture is also known when $G$ has domination number $1,2,$ or $3$ and $H$ is any graph. At this point, this is an exhaustive list of graph classes known to satisfy the conjecture. We note here that claw-free graphs are not known to be in any of these classes. These results and many others can be viewed in the survey \cite{VizingSurvey}

\medskip

This paper addresses a conjecture which is related to that of Vizing. In 2006, B.~Bre\v{s}ar, M. A. Henning, and S. Klav\v{z}ar \cite{BHK} posed a weaker version of Vizing's conjecture. Using the integer domination function of \cite{DHLF}, they defined integer domination, which we review later in this paper, and noticed that for any $G$, $\gamma(G)\le \gamma_{\{2\}}(G) \le \gamma_{\{3\}}(G) \le \dots$. In fact, the distances between these invariants are quite large. It is therefore surprising that the inequality $\gamma_{\{k\}}(G\square H)\ge \gamma(G)\gamma(H)$ has remained open for these last thirteen years for any $k$. The smaller the $k$, the closer to Vizing's conjecture, but so far, no one has proven this statement even for $k=1000$. The ``best'' value of $k$ in this weakened conjecture is $2$.

\medskip

Notice that if $G$ satisfies Vizing's conjecture for any $H$, then it also satisfies this integer domination conjecture. That is, if for a given $G$, $\gamma(G\square H)\ge \gamma(G) \gamma(H)$, since $\gamma_{\{2\}}(G\square H)\ge \gamma(G\square H)$, Bre\v{s}ar, Henning, and Klav\v{z}ar's integer domination conjecture is true for that $G$. This means that the conjecture is known for all classes of graphs from our list above. 

We now state the questions formally and note here that a positive answer to Question $1$ implies a positive answer to Question $2$. These questions also appear in \cite{VizingSurvey}.

\begin{quest}{\rm (Bre\v{s}ar, Henning, Klav\v{z}ar \cite{BHK})}\label{q1}
For any graphs $G$ and $H$, is it true that
\[\gamma_{\{2\}}(G\cart H) \geq \gamma(G)\gamma(H)?\]
\end{quest}

\begin{quest}{\rm (Bre\v{s}ar, Henning, Klav\v{z}ar \cite{BHK})}\label{q2}
Is there a natural number k such that for any pair of graphs $G, H$,
\[\gamma_{\{k\}}(G\cart H) \geq \gamma(G)\gamma(H)?\]
\end{quest}

In this note we answer these questions in the affirmative for claw-free $G$ and arbitrary $H$ by a vertex labeling method first introduced in~\cite{Krop}. We also show that for claw-free $G$ and any $H$, $\gamma_2(G\cart H)\ge \gamma(G)\gamma(H)$, where $\gamma_2(G)$ is the $2$-domination number of $G$, which we define in the next subsection. Our result for claw-free graphs is the first in fourteen years which shows the conjecture true for a class of graphs not known to satisfy Vizing's conjecture.

\noindent\textbf{Definitions and Notation.}
In this note, all graphs will be considered finite  and simple. Specifically, let $G$ be a graph with vertex set $V(G)$ and edge set $E(G)$. The order and size of $G$ will be denoted by $n(G) = |V(G)|$ and $m(G) = |E(G)|$, respectively. Two vertices $v$ and $w$ in $G$ are adjacent, or neighbors, if $vw\in E(G)$. A set of pairwise non-adjacent vertices in $G$ is an \emph{independent set}, or \emph{stable set}. The open neighborhood of a vertex $v\in V(G)$, written $N_G(v)$, is the set of all neighbors of $v$, whereas the closed neighborhood of $v$ is $N_G[v] = N_G(v)\cup \{v\}$. Let $S\subseteq V(G)$ and $v\in S$. The \emph{open} $S$-\emph{private neighborhood} of $v$ is defined as $\pn(v,S) = \{w\in V(G): N_G(w)\cap S = \{v\}\}$. A graph $G$ is called \emph{claw-free} if $G$ contains no $K_{1,3}$ as an induced subgraph. 

As mentioned previously, a set of vertices $S\subseteq V(G)$ is dominating if every vertex not in $S$ has a neighbor in $S$. If $S$ is a dominating set with the additional property that $S$ is also an independent set, then $S$ is a \emph{independent dominating set}. The cardinality of a minimum independent dominating set in $G$ is the \emph{independent domination number} of $G$, denoted $i(G)$. 

A set of vertices $S$ in a graph $G$ is \emph{2-dominating} if every vertex not in $S$ has at least 2 neighbors in $S$. The \emph{2-domination number} of $G$, written $\gamma_2(G)$, is the cardinality of a minimum 2-dominating set in $G$. With this definition, it is clear that every 2-dominating set is also a dominating set, and so, $\gamma_2(G) \ge \gamma(G)$.

For any integer $k\ge 1$, a function $f : V (G) \rightarrow \{0, 1, . . ., k\}$ is called a \emph{$\{k\}$-dominating function} if the sum of its function values over any closed neighborhood is at least $k$. The weight of a $\{k\}$-dominating function is the sum of its values over all the vertices. The $\{k\}$-domination number of $G$, $\gamma_{\{k\}}(G)$, is defined to be the minimum weight taken over all $\{k\}$-domination functions. 

The following invariant was introduced in \cite{CHHM} for $k=2$ and is also known as \emph{Italian domination number} in \cite{HK}.

\begin{defn}
For any integer $k\ge 1$, a function $f : V (G) \rightarrow \{0, 1, . . ., k\}$ is called a \emph{Roman $\{k\}$-dominating function} if for $v\in V(G)$ such that $f(v)= 0$, the sum of the function values over the closed neighborhood of $v$ is at least $k$. The weight of a Roman $\{k\}$-dominating function is the sum of its values over all the vertices. The Roman $\{k\}$-domination number of $G$, $\gamma_{Rk}(G)$, is defined to be the minimum weight taken over all Roman $\{k\}$-dominating functions.
\end{defn}

We note that we can make a formulation of a Roman $\{k\}$-dominating function on $G$ by replacing any vertex $v$ with $f(v)=\ell$, for $1\le \ell\le k$, by a clique of size $\ell$ where every vertex in the clique has function value $1$. That is, for any vertex $v$ so that $f(v)=\ell>0$, we replace $v$ by $v_*^{\ell}=K_{\ell}$ where every vertex $x$ in the clique $v_*^{\ell}$ has the same neighbors as $v$ and $f(x)=1$. We call the vertices $x$, \emph{copies} of $v$, and the graph produced by such a replacement of vertices, $G^*$. We call vertices $u\in V(G^*)$ so that $f(u)>0$, \emph{dominating vertices of $G^*$}. Furthermore, we call such a distribution of dominating vertices for a Roman $\{k\}$-dominating function, a \emph{Roman \{k\}-dominating set} of $G$. If $f$ is a Roman $\{k\}$-domination function of minimum weight, we call the distribution of dominating vertices a \emph{minimum Roman $\{k\}$-dominating set}.

The \emph{Cartesian product} of two graphs $G(V_1,E_1)$ and $H(V_2,E_2)$, denoted by $G \cart H$, is a graph with vertex set $V_1 \times V_2$ and edge set $E(G \cart H) = \{((u_1,v_1),(u_2,v_2)) : v_1=v_2 \mbox{ and } (u_1,u_2) \in E_1, \mbox{ or } u_1 = u_2 \mbox{ and } (v_1,v_2) \in E_2\}$.

For any graphs $G$ and $H$, let $f$ be a Roman $\{k\}$-dominating function of $G\cart H$. For any such $f$, we now extend the Cartesian product to the \emph{Roman Cartesian product} in which every vertex $v\in V(G\cart H)$, so that $f(v)=\ell>0$, is replaced by the clique $v_*^{\ell}$ as described previously. We denote this distinction by the notation $\cartr$ for the Roman Cartesian product. For any vertex $v\in V(G\cart H)$, the \emph{multiplicity} of $v$ is $\ell$ if $v$ was replaced by $v_*^{\ell}$ to form $V(G\cartr H)$.

Given a Roman $\{2\}$-dominating set $D$ of of $G\cartr H$, let $S$ be a subset of $D$ so that\\
\noindent $S=\{(g_1,h_1),\dots, (g_n,h_n)\}$ for $g_i\in V(G)$ and $h_i\in V(H)$, $i\in [|S|]$. We define the \emph{projection of $S$ onto $H$} as the set of vertices in $H$, $\{h_1,\dots, h_n\}$. Likewise, the \emph{projection of $S$ onto $G$} is the set of vertices in $G$, $\{g_1,\dots, g_n\}$.

For other graph theoretic terminology and definitions, we will typically follow~\cite{MHAYbookTD}.

For a given positive integer $k$, we will also make use of the standard notation $[k] = \{1,\dots, k\}$. 

\section{Main Result}
In this section we will prove our main result, but before doing so we will need a useful theorem of Allan and Laskar~\cite{Allan} that states equivalence between the domination number and independent domination number of claw-free graphs. We state this theorem formally as follows. 
\begin{thm}{\rm (Allan and Laskar \cite{Allan})}\label{claw-free thm}
If $G$ is a claw-free graph, then $\gamma(G) = i(G)$. 
\end{thm}

We will also make use of the following observation.
\begin{ob}\label{claw-free ob}
If $G$ is a claw-free graph and $S$ is a minimum independent dominating  set in $G$, then any vertex not in $S$ is adjacent to either one or two vertices in $S$. 
\end{ob}

We are now ready to present our main result. 
\begin{thm}
If $G$ is a claw-free graph and $H$ is a graph, then
\begin{equation*}
\gamma_{R2}(G\cart H)\ge \gamma(G)\gamma(H).
\end{equation*}
\end{thm}

\proof Let $G$ be a claw-free graph and $H$ be a graph. By Theorem~\ref{claw-free thm}, we may choose an independent dominating set $S\subseteq V(G)$ with $|S| = \gamma(G)$. For notational simplicity, suppose $|S| = k$, and label the vertices of $S$ by $v_1,\dots v_k$. Let $f$ be a Roman $\{2\}$-dominating function of $G\cart H$ and for this $f$ we consider the Roman Cartesian product, $G\cartr H$. Let $D\subseteq V(G\cartr H)$ be a minimum Roman $\{2\}$-dominating set of $G\cartr H$, and so, $|D| = \gamma_{R2}(G\cart H)$. 

We next devise a labeling scheme for the vertices in $D$ which is split into two separate parts; an initial labeling and a finishing labeling. The following labeling is the initial labeling, and we note that by Observation~\ref{claw-free ob} the following labeling will assign at most 2 entries to each label. For any vertex $u\in D$, we may refer to $u$ by its coordinates in the Cartesian product. That is, if $u$ is the vertex formed by the Cartesian product of $v\in V(G)$ and $h\in V(H)$, then we may write $u=(v,h)$.

\begin{enumerate}
\item[(1)] For $i\in [k]$, if $(v,h)\in D$ with $v\in \{v_i\}\cup \pn(v_i,S)$, then label $(v,h)$ by $\{i\}$. 

\item[(2)] For distinct $i,j\in [k]$, if $(v,h)\in D$ is a vertex with multiplicity $1$ with $v$ adjacent to both $v_i$ and $v_j$ in $G$, where $(v_i, h') \notin D$ for any $h'\in N_H[h]$, and $(v_j, h'')\notin D$ for any $h''\in N_H[h]$, then label $(v,h)$ by $\{i,j\}$. If $(v,h)$ is a vertex of $D$ with multiplicity $2$ with the same conditions, then label one copy by $\{i\}$ and the other by $\{j\}$.

\item[(3)] For distinct $i,j\in [k]$, if $(v,h)\in D$ with $v$ adjacent to both $v_i$ and $v_j$ in $G$, where $(v_i, h') \notin D$ for any $h'\in N_H[h]$, and $(v_j, h'')\in D$ for some $h''\in N_H[h]$, then label $(v,h)$ by $\{i\}$. 

\item[(4)] For distinct $i,j\in [k]$, if $(v,h)\in D$ with $v$ adjacent to both $v_i$ and $v_j$ in $G$, where $(v_i, h') \in D$ for some $h'\in N_H[h]$, and $(v_j, h'')\in D$ for some $h''\in N_H[h]$, then label $(v,h)$ by $\{i\}$ or $\{j\}$ arbitrarily. 
\end{enumerate}

If $S$ is an independent dominating set of $G$ so that all vertices of $G$ are dominated by exactly once vertex of $S$ (also known as a perfect dominating set), then all vertices in $G\cartr H$ are assigned labels according to (1) in the above process. In this case the final labeling scheme is not necessary and we may bypass it. The following is the finishing labeling. 
 \begin{enumerate}
 \item[(5)] If $(u,h)$ and $(v,h')$ are vertices in $D$ that are both labeled $\{i,j\}$ for distinct $i,j\in [k]$, and $hh'\in E(H)$, then relabel $(u,h)$ by $\{i\}$, and $(v,h')$ by $\{j\}$. 
 
\item[(6)] If $(u,h)$ and $(v,h')$ are vertices in $D$ that are labeled $\{i\}$ and $\{i,j\}$, respectively, for distinct $i,j\in [k]$, and $hh'\in E(H)$, then relabel $(v,h')$ by $\{j\}$. 

\item[(7)] If $(u,h)$ and $(v,h)$ are vertices in $D$ that are both labeled $\{i,j\}$ for distinct $i,j\in [k]$, then relabel $(u,h)$ by $\{i\}$, and $(v,h)$ by $\{j\}$.

\item[(8)] If $(u,h)$ and $(v,h)$ are vertices in $D$ with labels $\{i\}$ and $\{i,j\}$, respectively, then relabel $(v,h)$ by $\{j\}$. 

\item[(9)] If $(u,h)$ and $(v,h)$ are vertices in $D$ with labels $\{i,j\}$ and $\{j,\ell\}$, respectively, then relabel $(v,h)$ by $\{\ell\}$. 

\item[(10)] If $(u,h)$ and $(v,h)$ are vertices of $D$ both labeled $\{i\}$, then we may relabel one of $(u,h)$ or $(v,h)$ by any other label.
\end{enumerate}

\begin{claim}\label{claim}
We may apply labelings (1) - (10) to $D$ and produce a labeling such that each vertex has a label with exactly one entry. 
\end{claim}
\proof Suppose that after applying labelings (1) - (10), there exists a vertex of $D$, say $(x_{i_1},h)$, which has been assigned the labeling $\{i_1,i_2\}$. Notice that vertices of $D$ may receive labels with two entries only in labeling (2).

According to the conditions of labeling (2), $(x_{i_1},h)$ has multiplicity $1$, $(v_{i_1},h')\cap D = \emptyset$ and $(v_{i_2},h')\cap D = \emptyset$ for all vertices $h'\in N_H[h]$. Since $D$ is a Roman $\{2\}$-dominating set of $G\cartr H$, all vertices not in $D$ have at least two neighbors in $D$, it follows that both $(v_{i_1},h)$ and $(v_{i_2},h)$ have at least one other neighbor distinct from $(x_{i_1},h)$ in $D$.

Suppose a vertex $(x_{i_2},h)\in D$ for some $x_{i_2}\in V(G)$, different from $(x_{i_1},h)$, is the neighbor of $(v_{i_2},h)$. First suppose that $(x_{i_2},h)$ is assigned a label that contains $\{i_2\}$ (or $\{i_1\}$). By labeling $(7)$ and $(8)$, two vertices such as $(x_{i_2},h)$ and $(x_{i_1},h)$ would receive labels with one entry, which contradicts the possibility of the label $\{i_1,i_2\}$ on $(x_{i_1},h)$.

Finally, we suppose that $(x_{i_2},h)$ has been assigned the label $i_3$, which is distinct from either $i_1$ or $i_2$, for some $i_3\in[k]$. Let $n$ be the minimal index so that for $2\leq \ell\leq n$, $(x_{i_{\ell-1}},h)$ and $(x_{i_{\ell}},h)$ are adjacent to $(v_{i_{\ell-1}},h)$, and $(x_{i_n},h)$ is adjacent to some vertex $(v_{i_m},h)$ for some $m\in [n-1]$, where $(x_{i_\ell},h)$ is labeled $i_{\{\ell+1\}}$ for $2\leq \ell \leq n-1$.  

We consider the cycle $(v_{i_m},h), (x_{i_m},h), (v_{i_{m+1}},h),\dots, (x_{i_n},h)$. Vertex $(x_{i_n},h)$ may be labeled by $\{i_{n}\}$, $\{i_m\}$, or $\{i_m,i_{n}\}$. If the label on $(x_{i_n},h)$ contains $\{i_{n}\}$, then by labeling (10), we may relabel $(x_{i_{n-1}},h)$ by $\{i_{n-1}\}$, and continue relabeling vertex $(x_{i_\ell},h)$ by $\{i_{\ell}\}$ for $2\leq \ell \leq n-1$. However, by labeling (8), this means $(x_{i_1},h)$ could be labeled by $\{i_1\}$. If the label on $(x_{i_n},h)$ contains $\{i_{m}\}$, then by labeling (10), we may relabel $(x_{i_\ell},h)$ by $\{i_{\ell}\}$ for $2\leq \ell \leq m$. Again, by labeling (8), this means $(x_{i_1},h)$ could be labeled by $\{i_1\}$.

Thus, $(x_{i_1},h)$ could be relabeled by a label with one entry.
\smallqed 

According to Claim~\ref{claim}, each vertex of $D$ has been assigned a label with a single entry. Choose $i\in[k]$, project all vertices of $D$ labeled $i$ onto $H$, and call the projected vertices $U=\{u_1, \dots, u_{\ell}\}$. 

\begin{claim}
The set $U$ is a dominating set of $H$.
\end{claim}
\proof By way of contradiction, suppose $U$ is not a dominating set of $H$; that is, there exists $h\in V(H)$ such that $h\notin N[U]$. This means that in $G\cartr H$, $(v_i,h)$ is dominated by some vertex $(v,h)$ of $D$ labeled by $\{j\}$ for some $j\in [k]$ with $j\neq i$. Labelings (3) and (4) could not have been applied in this case, since $h\notin N[U]$. If labeling (2) had been applied to $(v,h)$, and then any or none of the labelings (5), (6), (7), (8), (9), or (10) had been applied, then $(v,h)$ would have been adjacent to some vertex $(u,h)$ labeled $\{i\}$. Since $(v,h)$ is adjacent to $(v_i,h)$ and $(v_j,h)$, labeling (1) does not apply. This produces a contradiction since no labeling could have been applied to $(v,h)$ but all vertices of $D$ are labeled by Claim 1. Thus, $U$ is a dominating set of $H$, and the proof of the claim is finished. \smallqed

By Claim 2, for each $i\in[k]$, we may project $D$ onto $H$ and obtain a dominating set of $H$. Recalling $\gamma_{R2}(G\cart H) = |D|$, $\gamma(G) = k$, and $|U| \ge \gamma(H)$, we observe the following,
\[
\begin{array}{lcl}
\gamma_{R2}(G\cart H) & = & \displaystyle{|D|} \1 \\
& \ge & \displaystyle{ \sum_{i=1}^{k}\gamma(H)} \1 \\
& = & \displaystyle{\gamma(G)\gamma(H)}.
\end{array}
\]
Thus, the proof of the theorem is complete. \qed

Since a $2$-dominating set produces a Roman $\{2\}$-dominating function, for any graph $G$, $\gamma_{R2}(G)\le \gamma_2(G)$. Furthermore, $\gamma_{R2}(G)\le \gamma_{\{2\}}(G)$ since the $\{2\}$-dominating function is a minimum over a set with more restrictions than the Roman $\{2\}$-dominating function. These two observations lead to the following.

\begin{cor}
For any claw-free graph $G$ and any graph $H$,
\begin{align*}
\gamma_{\{2\}}(G\cart H)\ge \gamma(G)\gamma(H)\\
\gamma_2(G\cart H)\ge \gamma(G)\gamma(H).
\end{align*}
\end{cor}

Concluding Remarks:

Since any graph $G$ which is known to satisfy Conjecture \ref{Vizing} for all graphs $H$ also answers Question \ref{q1} in the affirmative, our future goal is to answer Question \ref{q1} for any other class of graphs which is not known to satisfy Vizing's conjecture. We see such results as steps forward to understanding the complications in solving Vizing's conjecture.

\end{document}